\documentclass[12pt]{article} 

\newcommand{\kmcomment}[1]{}
\usepackage[leqno]{amsmath}
\usepackage{amscd}
\usepackage{scalerel}
\usepackage{newtxtext}

\usepackage{mathtools}
  \newcommand{\ds}{\ensuremath{\displaystyle }}
  \newcommand{\pdel}{\partial} 
  \newcommand{\myHom}[1]{\textrm{H}_{#1}}
  \newcommand{\mR}{\ensuremath{\mathbb{R}}} 
  \newcommand{\mZ}{\ensuremath{\mathbb{Z}}} 
  \newcommand{\xb}[1]{x_{#1}}
  \newcommand{\frakg}{\mathfrak{g}}
   \newcommand{\tbdl}[1]{\mathrm{T}(#1)}
  \newcommand{\cbdl}[1]{\mathrm{T}^{*}(#1)}
  \newcommand{\Sbt}[2]{[#1,#2]}
  \newcommand{\SbtS}[2]{[#1,#2]_{\scriptstyle \text{S}}} 
  \newcommand{\SbtE}[3]{\ensuremath{[#2,#3]_{#1}}} 
  \newcommand{\parity}[1]{(-1)^{#1}}
  \newcommand{\CSp}[1]{\text{C}_{#1}} 
  \newcommand{\mw}{\mywedge}
  \newcommand{\we}{\wedge}
  \newcommand{\wfg}[1]{ \frakg_{<#1>}} 
  \newcommand{\myCSW}[2]{\text{C}_{#1,#2}} 
  \newcommand{\mywedge}{\bigtriangleup} 
  \newcommand{\myCS}[1]{\text{C}_{#1}} 
\newcommand{\Eu}{E}
\newcommand{\SW}[1]{SW(#1)\;}
\newcommand{\PW}[1]{PW(#1)\;}

\newcommand{\kmMat}[2]{{\scriptscriptstyle \left[\begin{array}{c} g( #1 )\\ #2 \end{array}\right]}}

\newcommand{\kmMatt}[2]{{\scriptstyle #1 }}

\newcommand{\I}{II}
\newcommand{\ActOne}{\mathcal{A}_{1}}
\newcommand{\partt}{\text{part}}  
\newcommand{\AAs}{As} \newcommand{\BBs}{Bs} \newcommand{\CCs}{Cs}

\usepackage{wrapfig2}

\usepackage{fancyvrb,shortvrb}
 \RequirePackage{fancyvrb}
  \DefineVerbatimEnvironment%
      {newverb}{Verbatim}
      {xleftmargin=10mm}
\DefineVerbatimEnvironment%
     {newverb*}{Verbatim}
     {xleftmargin=10mm,showspaces=true,showtabs=true}

\usepackage{listings}
\RequirePackage{listings}
\lstnewenvironment{kmverb}
     { \lstset{xleftmargin=10mm, breaklines=true}}
     { }
\lstnewenvironment{kmverbb}[1]
     { \lstset{xleftmargin=#1, showspaces=false, breaklines=true}}
     { }
 
\usepackage[a4paper]{geometry}
\usepackage{amsthm}
\theoremstyle{definition} 
\newtheorem{theorem}{Theorem}
\newtheorem*{theorem*}{Theorem}
\newtheorem{definition}[theorem]{Definition}
\newtheorem*{definition*}{Definition}

\newtheorem{prop}{Proposition}[section]
\newtheorem{remark}{Remark}[section]
\newtheorem{lemma}{Lemma}
\newtheorem{Cor}{Corollary}

\title{Super homology groups of differential forms and
    vector fields on Euclidean line}

\author{Kentaro Mikami \and Tadayoshi Mizutani}

\usepackage{setspace}
\numberwithin{equation}{section}

\allowdisplaybreaks

\newcommand{\almostprimary}{}
\newcommand{\ddx}{\frac{d}{dx}}

 \usepackage{newtxmath}
\AtBeginDocument{\parindent=0pt \setlength{\headheight}{21pt}}
\begin{document}
\maketitle

\begin{abstract}
The arXiv:2105.09738 claims several stuffs. In particular, we recall the
following two.

(1) Vector fields and differential forms become a Lie superalgebra structure
for each manifold.  

(2) For an n-dimensional Euclidean space,  
vector fields and differential forms with polynomial coefficients become a 
double weighted Lie superalgebra.  By using Euler vector field, 
the Betti numbers are 0 except the last one 
if the primary weight and the secondary weight are different.   

Now, a simple question arises: What happens when  
the primary weight and the secondary weight are equal?  
This note shall give a complete answer to the question for the case $n=1$.  

MSC(Math Subject Classification)2020  are 17B10, 17B15.  
\end{abstract}

\kmcomment{ 
17Bxx: (Lie algebras and Lie superalgebras)
17B01: (Identities, free Lie (super)algebras)
17B05: (Structure theory for Lie algebras and superalgebras)
17B08: (Coadjoint orbits; nilpotent varieties)
17B10: (Representations of Lie algebras and Lie superalgebras, algebraic theory (weights))
17B15: (Representations of Lie algebras and Lie superalgebras, analytic theory)
}


\section{Introduction}
Let \(M\) be a differential manifold of dimension \(n\). 

A typical Lie superalgebra structure on 
\( \tbdl{M} \oplus \Lambda^{0} \cbdl{M} \oplus \cbdl{M} \oplus \cdots 
\oplus \Lambda^{n} \cbdl{M} \) 
is defined by 
\[ 
    \Sbt{X}{Y} = \SbtS{X}{Y} \quad \text{(Schouten bracket)}, 
    \quad 
    \Sbt{X}{\alpha} = L_{X} \alpha \quad \text{(Lie derivative)},  
    \quad 
    \Sbt{\alpha}{\beta} = \parity{a} d( \alpha \we \beta ) \; , \]
    where \(X,Y\) are vector fields, \(\alpha, \beta\) are differential
    forms and \(\alpha \) is \(a\)-form. 

Assume now \( M = \mR^n \) with Cartesian coordinate \( x_1, \ldots, x_n\).

Then we may restrict vector fields spanned by \( \frac{\pdel}{ \pdel x_1 },
\ldots,  \frac{\pdel}{ \pdel x_n } \) whose coefficients, which are only
\((j+1)\)-homogeneous  polynomials of \(x_1, \ldots, x_n\), which we denote
\( \wfg{0,j}\).  Also,  we restrict differential \(i\)-forms which are
generated by \( d x_1, \ldots, d x_n \) with coefficients of 
\((j+1)\)-homogeneous polynomials of \(x_1, \ldots, x_n\) too, 
which is denoted by \( \wfg{-i-1,j}\).  


It is known that those objects are invariant under our super bracket as
follows. 
\begin{prop}
\(\ds \SbtE{}{ \wfg{i,h} } { \wfg{j,k} } \subset  
 \wfg{i+j,h+k}  
\) holds for \( i,j \leqq 0\) and \( h,k \geqq -1\).  
\end{prop}
Thus, we have doubly weighted chain complex 
\(\myCSW{\bullet}{w,h} := \sum_{m=mL}^{mH} \myCSW{m}{w,h} \),  where  
\[\myCSW{m}{w,h}= \prod_{i=1}^{m} 
    \{\wfg{ w_{i}, h_{i} } \mid \sum_{i} w_{i}
= w \text{ (the primary weight), and } \sum_{i} h_{i} = h 
\text{ (the secondary weight)} \}
\;.
\]  

\begin{theorem}
The Betti numbers except the last are 0 when \( w-h \ne 0\) for 
the double weighted homology groups
of the above double weighted chain complex \(\myCSW{m}{w,h}\)
of \(\mR^{n}\). The last Betti number is given by the alternating sum of
chain space dimensions.  (Cf.\ \cite{Mik:Miz:superForms})
\end{theorem}
\textbf{Proof:}
Using the Euler vector field 
\(\ds \Eu = \sum_{i=1}^{n} \xb{i} \frac{\pdel}{\pdel \xb{i}} \), we have 
\begin{align*}
\Sbt{\Eu} {x^{P} d x^{A}} &= (|P|+ |A|) x^{P} d x^{A} = 
 (\SW{x^{P} d x^{A}}- \PW{x^{P} d x^{A}} ) x^{P} d x^{A}  \\
\Sbt{\Eu} {x^{P} \pdel/ \pdel  \xb{j}} &
= (|P|- 1) x^{P}\pdel /\pdel\xb{j}
= \SW{x^{P}\pdel /\pdel\xb{j} } x^{P}\pdel /\pdel\xb{j}
\end{align*}
where 
\(
 \PW{\alpha} \) and  
\(\SW{\alpha}\) mean the primary and the secondary weight of \( \alpha \).  

Since the double weight of \(\Eu\) is \((0,0)\), 
the definition of the boundary operator implies 
\begin{align*}
\pdel ( \Eu \mw Y_{1} \mw \cdots \mw Y_{m} ) &= - 
\Eu \mw \pdel( Y_{1}  \mw \cdots \mw Y_{m} ) + \sum_{i=1}^{m} 
Y_{1} \mw \cdots \mw \Sbt{\Eu}{Y_{i}}\mw\cdot \mw Y_{m}  
\\\shortintertext{for \( mL \leqq m < mH\). 
When  \( Y_{i} \in \wfg{w_{i},h_{i}}\), 
    \(\Sbt{\Eu}{Y_{i}} = (h_{i}-w_{i} ) Y_{i}\) 
because \(\Eu\) is the Euler vector field, so we have   
}
\pdel ( \Eu \mw Y_{1} \mw \cdots \mw Y_{m} ) &= - 
\Eu \mw \pdel(Y_{1}\mw \cdots \mw Y_{m} ) + \sum_{i=1}^{m} ( h_{i} - w_{i} ) 
Y_{1} \mw \cdots \mw Y_{m}  
\\ &= - 
\Eu \mw \pdel( Y_{1}  \mw \cdots \mw Y_{m} ) + (h-w)
Y_{1} \mw \cdots \mw Y_{m} 
\end{align*}
for \( Y_{1} \mw \cdots \mw Y_{m} \in \myCSW{m}{w,h}\). 
Thus, if \( h-w \ne 0\), we see each \(m\)-cycle is exact and so the Betti
number is 0. The last Betti number is given by
\( (-1)^{mH} \sum_{m=mL}^{mH} (-1)^ m \dim  \myCSW{m}{w,h}\). 
\qed

\bigskip



From now on, we introduce a new notation for 
 \(m\)-th chain space as follows.   
\begin{definition}
\( \kmMat{0}{j}  :=  \wfg{0, j-1}\), which is 
the space of vector fields with \(j\)-th homogeneous
polynomials. 

\( \kmMat{i}{j}  :=  \wfg{-i, j-1}\) for \(i=1, \ldots, n+1\) and \( j \geq
0\), which is the space of \((i-1)\)-differential forms 
with \(j\)-th homogeneous polynomials.  

For a given \(m\)-th chain 
\begin{equation}
 \kmMat{i_1}{j_1} \mw \kmMat{i_2}{j_2} \mw \cdots\mw  \kmMat{i_m}{j_m}
 \quad \text{where}  \quad 
 i_{s} \in \{ 0,1,\dots, n+1\} \quad\text{and}\quad j_{s} \geq 0 \;, 
 \label{second:expression}
\end{equation}
the primary and secondary weights 
are defined by  
\( w = - \sum_{\ell=1}^{m} {i_\ell}\) 
and \( h  = -m + \sum_{\ell=1}^{m} j_{\ell}\).   
In this note, we redefine the primary weight by 
\( \sum_{\ell=1}^{m} {i_\ell}\), namely only change the sign simply.   

By the super symmetric property of \(\mw\), i.e., 
\( A \mw A' + \parity{i i' } A' \mw A = 0 \) holds for 
\( A\in \kmMat{i}{j} \) and \( A'\in \kmMat{i'}{j'} \),  
we reorder 
\(i_1, i_2, \ldots, i_m\) non-decreasing order, namely 
\(i_1 \leq  i_2 \leq \cdots \leq  \ldots \leq i_m\)
in \eqref{second:expression}.   
Suppose 
\(i_1 =  i_2 =  \ldots = i_p\) for instance. Then we reorder their second
indices so that    
\(j_1 \leq  j_2 \leq \cdots \leq  \ldots \leq j_p\).   
We use the natural ordering: 
first compare by type of forms or vector fields, then compare by the
polynomials.   
In particular, when \(n=1\) case, if 
\(i_1 =  i_2 =  \ldots = i_p\) and  it is a common even integer, then 
\( j_1, \ldots, j_p\) should be mutually distinct, namely 
\(j_1 <  j_2 < \cdots < j_p\).  
\end{definition}

\section{Research on the case where the spatial dimension is one}
Here, we restrict the dimension of the space we handle now is 1.  
\kmcomment{ 
Although
in the previous arXiv paper, the primary weight of \(q\)-form is \(- (1 + q)
\), here the primary weight of forms are
positive,  namely, the primary weight of \(q\)-form is \( 1 + q \).

The space is \(\mR\) and we handle vector fields, polynomials, or
differential forms with homogeneous polynomials like as \( f(x) \frac{d}{d x}, g(x), h(x) d x\) where \( f(x), g(x), h(x) \) are 
homogeneous polynomial functions of \(x\). 

Let \(f\) be a \(m\)- homogeneous function and \(E = x \frac{d}{dx}\), 
which is called the Euler vector field. Then we get 
\begin{equation}
    \Sbt{E}{f} = m f\ , \quad 
    \Sbt{E}{f dx} = (m+1) f dx \ , \quad 
    \Sbt{E}{f \frac{ d}{dx}} = (m-1) f \frac{d}{dx}\ . 
    \end{equation}
} 

\kmcomment{
Case of \(w = - h = 0\): 
\(\myCSW{1}{w,h}=  \wfg{ 0 , 0 } \),  
\(\myCSW{2}{w,h}=  \wfg{ 0 , -1 } \mw  \wfg{ 0 , 1 } 
+  \wfg{ 0 , 0 } \mw  \wfg{ 0 , 0 } \) and because of 
the 2nd term is 0 imply  
\(\myCSW{2}{w,h}=  \wfg{ 0 , -1 } \mw  \wfg{ 0 , 1 } \).   

\[ \left[\normalsize \begin{array}{c} A  \\ B \end{array} \right] \quad 
 \left[\small \begin{array}{c} A  \\ B \end{array} \right] \quad  
 \left[\footnotesize \begin{array}{c} A  \\ B \end{array} \right] \quad  
 \left[\scriptsize \begin{array}{c} A  \\ B \end{array} \right] \quad 
 \left[\tiny \begin{array}{c} A  \\ B \end{array} \right] \quad 
\kmMat{A}{B} \quad 
\kmMat{A}{B}\]
}

\kmcomment{ 
The primary weight of  m-th  multi-vector field  is \(m-1\) and the primary
weight of \(\ell\)-th differential form should be \( -\ell -1\).  However,
when the multivector fields do not appear, even if we redefine the primary
weight of \(\ell\)-th differential form by \( \ell+1\) there is no
confusion.  Thus, here after we use the second definition of primary weight
of differential forms, i.e., those values are positive. 
}

We study the homology groups of chain complex \( \{ \myCSW{m}{w,h=-w} \}_{m}\) of the
superalgebra associated with \(\mR\), i.e., when \(n=1\). 
Then only \( g(0), g(1), g(2) \) appear and we define the primary weight of \(
g(i)\) is \(i\).


Since \( P \mw Q + Q \mw P =0\) for 
\( P = \kmMat{0}{j}  \text{ or } \kmMat{2}{k} \), 
and 
 \( P \mw Q = Q \mw P \) for 
\( P= \kmMat{1}{j}, Q=\kmMat{1}{k}\),  we re-order 
\( \kmMat{i_1}{j_1} \mw \kmMat{i_2}{j_2} \mw \cdots \mw \kmMat{i_m} {j_m} \)
to  
\begin{equation} \kmMat{0}{i_1} \mw \cdots\mw \kmMat{0}{i_a} \mw 
\kmMat{1}{j_1} \mw \cdots\mw \kmMat{1}{j_b} \mw 
\kmMat{2}{k_1} \mw \cdots\mw \kmMat{2}{k_c} 
\label{imadake:one}
\end{equation} 
where  \( i_1 < \cdots < i_a\),  \( j_1 \leqq  \cdots \leqq  j_b\), and 
\(k_1 < \cdots < k_c\).  

We recall a small lemma.
\begin{lemma}
Consider a finite well-ordered sequence of  distinct integers 
\( a_1 < a_2 < \cdots < a_m  \).  
Then their mean satisfies 
\begin{equation}
    a_{1} + \frac{ (m-1) }{2} \leqq \frac{1}{m}\sum_{i=1}^{m} a_i \leqq  a_{m} - 
     \frac{ (m-1) }{2} 
\end{equation}
\end{lemma}
Proof:  We will show the following by induction:    
\( \displaystyle 
    m a_{1} + \frac{m (m-1) }{2} \leqq \sum_{i=1}^{m} a_i \leqq m a_{m} - 
     \frac{ m (m-1) }{2} 
 \). 

Define a sequence \( L_i := a_1 + i-1\) for \(i=1,...,m\), which is monotone
increasing of slope ratio 1. 
\( L_1 = a_1\) holds and  we see that \( L_i \leq a_i \) for
\(i=1,\ldots, m \) by induction. Thus, \( \sum_{i=1}^m a_i \geqq 
\sum_{i=1} ^m L_i \)  and 
\( \sum_{i=1} ^m L_i = m a_1 + \frac{m(m-1)}{2} \).   

Define a sequence \( U_i := a_m -m + i\) for \(i=1,...,m\), which is monotone
increasing of slope ratio 1. By the same discussion above, we see that  
 \( a_i \leq U_i \) for \(i=1,\ldots, m \) by induction, and 
 \( \sum_{i=1}^m a_i \leqq 
\sum_{i=1} ^m U_i =  m a_m + \frac{m(m-1)}{2} \).   
\qed 

\kmcomment{
\[\Sbt{A}{B} \,  \SbtD{A}{B}, \TCkt{X}{Y} \] 
\( \SbtT{A}{B}, \Sbtt{A}{B},\SbtN{A}{B},\SbtSZ{A}{B},\SbtE{00}{A}{B} \)
\[ \SbtT{A}{B}, \Sbtt{A}{B}, \SbtN{A}{B},\SbtSZ{A}{B}, \SbtE{0}{A}{B} \]
}

\begin{theorem} \label{thm:n:1}
Assume \(n=1\) and \(w = -h \). The chain complex  
\(\{ \myCSW{m}{w,-w} \}_{m} \) has non-trivial chain spaces for \( m= w, w+1, w+2,
w+3\) for \(w > 0\). When \(w=0\),    
non-trivial chain spaces are given by \( m= w+1, w+2, w+3\).   
\end{theorem}
{\textbf 
Proof:} 
Let us recall \eqref{imadake:one}.  Then, we get 
\begin{align}
    m & := a+b+c \label{my:m}  \\
    w & := 0 \cdot a +1 \cdot  b +2 \cdot c = b+2 c \label{my:w} \\
    \shortintertext{
    about secondary weight \(h\), we see}
    M & := m + h = m - w = a-c \label{my:2nd} \quad \text{by the assumption }
    w+h=0  \\  
    M & \geq  (0+1+\cdots + a-1) + 0\cdot b + (0+1+\cdots+ c-1)= (a-1)a/2
    +(c-1)c/2\;. \label{my:M}
    \\\shortintertext{The above two relations imply }
      & (a-1)a/2 +(c-1)c/2 \leq a-c\;, \ \text{i.e.,} \  0 \leq c(c+1) \leq -
      a(a-3)\; . \label{my:N}
\end{align}
The above 
\eqref{my:N} shows \( a=0,1,2,3\). 
If \(a=0\) or 3 
then \(c=0\), else if \(a=1\)  or 2 then  \(c=0\) or  1 hold.   
\kmcomment{
Assume \(c=0\). Then \( a = 0,1,2,3 \) and \( m=a+b, w=b, M \geqq (a-1)a/2, M =m-w = a \). Then 
\( M \geqq (M-1)M/2\) and so \(0 \leqq M \leqq 3\). 
Assume \(a=0\)  as another approach. 
Then \( m=b+c, w=b+2c\) and \(c=w-m = - M, M \geqq (c-1)
c/2\). Then 
\( M \geqq (M+1)M/2\) and so \(0 \leqq M \leqq 1\). Thus, we conclude \( 0
\leqq M \leqq 3\), namely \( m=w+0, w+1, w+2, w+3\). In particular, when \(w=0\),
\( m= 1, 2, 3\) removing \(m=0\).  
}
Thus,
using the relation \( m-w=a-c\) in \eqref{my:2nd},  
we conclude \( m=w+0, w+1, w+2, w+3\). In particular, when \(w=0\),
\( m= 1, 2, 3\) removing \(m=0\).  
\qed 

We recall here the general definition of the boundary operator of general Lie
superalgebra. 
\begin{definition}
For a \(\mZ\)-graded Lie superalgebra 
\(\ds \frakg = \sum_{i} \frakg_{i} \), $m$-th chain space is defined by 
\(\CSp{m} = \underbrace{\frakg  \mw  \cdots  \mw  \frakg}_{m}\)
and the boundary operator is given by 
\begin{align} 
\pdel (A_{1}  \mw  A_{2}  \mw  
\cdots  \mw  A_{m}) 
=& \sum_{i<j} \parity{ i-1 + a_{i}\sum_{i<s<j} a_{s} }
A_{1}  \mw  \cdots 
\widehat{ A_{i} }\cdots  \mw  
\Sbt{A_{i}}{A_{j}}  \mw  
 \cdots  \mw  A_{m}
 \label{bdary:orig}
\end{align}
where \( A_{i} \in \frakg_{a_{i}}\). We see that \( \pdel \circ \pdel = 0\)
holds, and have the homology groups by 
\[
\myHom{m} (\frakg ) =  \ker(\pdel_{} : \myCS{m} \rightarrow
\myCS{m-1})/ \pdel ( \myCS{m+1} )\;.  
\]
\end{definition}

\subsection{\(n=1, w=-h=0\)} 
By Theorem \ref{thm:n:1}, \(m\) of the chain space \( \myCSW{m}{0,-0} \)  
 runs 1,2,3.  
\begin{align*}
    \myCSW{1}{0,-0} &=  \kmMat{0}{1} =\mR\left(  x \frac{d}{d x}\right) \\  
    \myCSW{2}{0,-0} &=  \kmMat{0}{0} \mw \kmMat{0}{2} = \mR \left(\frac{d}{d x} \mw x^2 \frac{d}{d x}\right) \\  
    \myCSW{3}{0,-0} &=  \kmMat{0}{0} \mw \kmMat{0}{1}\mw \kmMat{0}{2} = \mR
 \left( \frac{d}{d x} \mw x \frac{d}{d x}  \mw x^2 \frac{d}{d x}\right) .   
\end{align*}
\begin{align*}
\pdel 
\left(\frac{d}{d x} \mw x^2 \frac{d}{d x}\right) = &  
\Sbt{ \frac{d}{d x}}{ x^2 \frac{d}{d x}} = 2 x \frac{d}{d x} \quad
\text{and} \\  
\pdel \left(
\frac{d}{d x} \mw x \frac{d}{d x}  \mw x^2 \frac{d}{d x}\right) =& 
\Sbt{ \frac{d}{d x}}{ x \frac{d}{d x} }  \mw x^2 \frac{d}{d x} + 
x \frac{d}{d x}  \mw\Sbt{ \frac{d}{d x}}{ x^2 \frac{d}{d x} } 
- \frac{d}{d x} \mw \Sbt{ x \frac{d}{d x} }{ x^2 \frac{d}{d x}} 
\\
= & 
 \frac{d}{d x} \mw x^2 \frac{d}{d x} + 
x \frac{d}{d x}  \mw 2 x  \frac{d}{d x}  
- \frac{d}{d x} \mw ( 
 2 x^2  \frac{d}{d x} - 
 x^2 \frac{d}{d x} ) = 0 \; .   
\end{align*}
Thus, the Betti numbers are 0,0,1. 


\subsection{\(n=1, w=-h=1\)} 

By Theorem \ref{thm:n:1}, \(m\) of the chain space \( \myCSW{m}{1,-1} \)  
 runs 1,2,3,4.  
\begin{align*}
\myCSW{1}{1,-1} &=  \kmMat{1}{0} =\mR \\  \myCSW{2}{1,-1} & 
=\kmMat{0}{0}\mw\kmMat{1}{1} + \kmMat{0}{1}\mw\kmMat{1}{0} 
= \mR \left(\frac{d}{d x} \mw x\right)+\mR \left(x\frac{d}{d x}\mw 1\right) \\ 
\myCSW{3}{1,-1} &=\kmMat{0}{0} \mw \kmMat{0}{1}\mw \kmMat{1}{1}+\kmMat{0}{0}
\mw \kmMat{0}{2}\mw \kmMat{1}{0} \\& 
= \mR\left( \frac{d}{d x} \mw x \frac{d}{d x}  \mw x \right)
                 + 
  \mR \left( \frac{d}{d x} \mw x^2 \frac{d}{d x}  \mw 1 \right) \\
\myCSW{4}{1,-1}&=\kmMat{0}{0}\mw\kmMat{0}{1}\mw\kmMat{0}{2}\mw\kmMat{1}{0} 
= \mR\left(\frac{d}{d x}\mw x \frac{d}{d x}\mw x^2 \frac{d}{d x}\mw 1 \right)\ . 
\end{align*}
The super bracket for forms are defined by \( \Sbt{\alpha}{\beta} =
\parity{a} d( \alpha \we \beta) \) where \(a= \text{usual degree of the form
} \alpha\) and \( \Sbt{X} {\beta} = 
L_{X} \beta \),  where  
\(L_{X}  =  i_{X} \circ d + d \circ i_{X}   \), which is called by the 
Lie derivation by \(X\).   
\begin{align*}
    \pdel ( \frac{d}{d x} \mw  x ) = &\Sbt{ \frac{d}{d x}} { x }  = 1\: ,  \quad  \pdel( x\frac{d}{d x}\mw 1 ) =  \Sbt{ x\frac{d}{d x} }{ 1 } = 0 
\\
  \pdel( \frac{d}{d x} \mw x \frac{d}{d x}  \mw x )  
  =&   \Sbt{  \frac{d}{d x}}{ x \frac{d}{d x} }  \mw x  +  
  x \frac{d}{d x} \mw \Sbt{ \frac{d}{d x} }{  x }   
  -    \frac{d}{d x} \mw \Sbt{ x \frac{d}{d x} }{ x }  
  \\=& 
     \frac{d}{d x}   \mw x  +  x \frac{d}{d x} \mw 1 - \frac{d}{d x} \mw x   
     =   x \frac{d}{d x} \mw 1 
 \\ 
  \pdel( \frac{d}{d x} \mw x^2 \frac{d}{d x}  \mw 1 )
  =& 
  \Sbt{ \frac{d}{d x} }{ x^2 \frac{d}{d x} } \mw 1 
  +   x^2 \frac{d}{d x}  \mw \Sbt{ \frac{d}{d x} }{ 1} 
  -  \frac{d}{d x} \mw \Sbt{ x^2 \frac{d}{d x} }{ 1 }
  = 2   x \frac{d}{d x} \mw 1 
  \\
%
\pdel(\frac{d}{d x}\mw x \frac{d}{d x}\mw x^2 \frac{d}{d x}\mw 1 ) 
= &  \Sbt{\frac{d}{d x} }{ x \frac{d}{d x} }  \mw x^2 \frac{d}{d x}\mw 1 
+ x \frac{d}{dx} \mw \Sbt{\frac{d}{d x}}{ x^2 \frac{d}{d x}}\mw 1
+ x \frac{d}{d x}\mw x^2 \frac{d}{d x}\mw \Sbt{\frac{d}{d x}}{ 1 }
 \\  & 
- \frac{d}{d x}\mw \Sbt{ x \frac{d}{d x}}{  x^2 \frac{d}{d x} } \mw  1  
- \frac{d}{d x}\mw  x^2 \frac{d}{d x}\mw \Sbt{ x \frac{d}{d x} }{1}  
+ \frac{d}{d x}\mw x \frac{d}{d x}\mw \Sbt{ x^2 \frac{d}{d x}}{ 1} 
\\
= &  {  \frac{d}{d x} }  \mw x^2 \frac{d}{d x}\mw 1 
+ x \frac{d}{d x}\mw {2 x \frac{d}{d x}}\mw 1    
+ x \frac{d}{d x}\mw x^2 \frac{d}{d x}\mw  0 \\  & 
- \frac{d}{d x}\mw {  x^2 \frac{d}{d x} } \mw  1  
- \frac{d}{d x}\mw  x^2 \frac{d}{d x}\mw 0
+ \frac{d}{d x}\mw x \frac{d}{d x}\mw 0 
\\
= &  { \frac{d}{d x} }  \mw x^2 \frac{d}{d x}\mw 1 
 + 0  
 + 0 
+  0  
- \frac{d}{d x}\mw {  x^2 \frac{d}{d x} } \mw  1  
-  0
+  0 = 
0\; . 
\end{align*}
Thus, the Betti numbers are 0,0,1,1. 

\subsection{\(n=1, w=-h=2\)} 
By Theorem \ref{thm:n:1}, \(m\) of the chain space \( \myCSW{m}{2,-2} \)
runs 2,3,4,5.  
\kmcomment{
The super bracket for forms are defined by \( \Sbt{\alpha}{\beta} =
\parity{a} d( \alpha \we \beta) \) and \( \Sbt{X} {\beta} = L_{X} \beta = 
( i_{X} \circ d + d \circ i_{X} ) \beta \). 
}
When \(m=2\), \(i_1+i_2=w=2\) and \(k_1 + k_2 = 0\) hold, and 
\( \myCSW{2}{2,-2}=\kmMat{0}{0} \mw  \kmMat{2}{0} + 
 \kmMat{1}{0} \mw  \kmMat{1}{0} \).  
We omit the detail explanation because the discussion is almost same, we
show the all of results only. 
\begin{align*}
\myCSW{2}{2,-2}=& \kmMat{0}{0} \mw  \kmMat{2}{0} + 
 \kmMat{1}{0} \mw  \kmMat{1}{0}\; , \\  
\myCSW{3}{2,-2}=& \kmMat{0}{0}\mw\kmMat{0}{1}\mw \kmMat{2}{0}
+ \kmMat{0}{0}\mw\kmMat{1}{0}\mw \kmMat{1}{1}
+ \kmMat{0}{1}\mw\kmMat{1}{0}\mw \kmMat{1}{0}
\\
\myCSW{4}{2,-2}=&\kmMat{0}{0}\mw\kmMat{0}{1}\mw \kmMat{1}{0}\mw \kmMat{1}{1}
+ \kmMat{0}{0}\mw\kmMat{0}{2}\mw \kmMat{1}{0}\mw \kmMat{1}{0}
\\
\myCSW{5}{2,-2}=&\kmMat{0}{0}\mw\kmMat{0}{1}\mw \kmMat{0}{2}\mw \kmMat{1}{0}
\mw  \kmMat{1}{0}\, .  
\end{align*}
In order to get the Betti numbers, we calculate each rank of the boundary
operator \(\pdel\). 
The rank on \( \myCSW{2}{2,-2} \) is 0, because the dead end of the boundary
 operator.   

About the rank on \( \myCSW{3}{2,-2} \), 
the boundary image of the first term vanishes as follows.
\begin{align*}
   \pdel (  \ddx \mw x\ddx \mw d x  )  
    &= \Sbt{ \frac{d}{dx} }{ x \frac{d}{dx}} \mw d x   
    + { x \frac{d}{dx}} \mw \Sbt{ \frac{d}{d x}}{d x}  - \frac{d}{dx}  \mw 
    \Sbt{ x \frac{d}{dx}}{ d x } \\ &= \frac{d}{dx} \mw (dx) + 0 -  \frac{d}{dx}
    \mw (dx) = 0 \; .  
\end{align*}
The boundary image of the 2nd term 
\( \kmMat{0}{0}\mw\kmMat{1}{0}\mw \kmMat{1}{1} \) is included in  
\( \myCSW{2}{2,-2}\) as follows.   
\begin{align*}
   \pdel (  \ddx \mw 1 \mw  x  )  
    &= \Sbt{ \frac{d}{dx} }{ 1 } \mw  x   
    + { 1 } \mw \Sbt{ \ddx }{ x}  - \ddx   \mw 
    \Sbt{ 1 }{ x } = 1 \mw 1  -  \ddx \mw d x
     \; ,  
\end{align*}
and the boundary image of the 3rd term is 0.   
Thus the rank on \( \myCSW{3}{2,-2} \)  is 1.  

About the rank on \( \myCSW{4}{2,-2} \): 
\begin{align*} 
    \pdel ( \ddx \mw x \ddx \mw 1 \mw x ) =& 
    \Sbt{\ddx}{ x \ddx } \mw 1 \mw x  
    +  x \ddx \mw \Sbt{\ddx}{ 1 }  \mw x  
    +  x \ddx \mw 1 \mw \Sbt{\ddx}{ x }   
                                        \\& 
-  \ddx  \mw ( \Sbt{x \ddx}{1} \mw x + 1 \mw \Sbt{x\ddx}{x} ) 
    + \ddx \mw x \ddx \mw \Sbt{1}{x}\\ 
    =&  x \ddx \mw ( 1 \mw 1 - \ddx \mw d x )
    =  x \ddx \mw  1 \mw 1 +  \ddx \mw x \ddx \mw d x 
    \\
    \pdel ( \ddx \mw x^2 \ddx \mw 1 \mw 1 ) =& 
    \Sbt{\ddx}{ x^2 \ddx } \mw 1 \mw 1 + 0 + 0   
    -  \ddx \mw ( 0+ 0) + \ddx \mw x^2 \mw 0  = 2 x \ddx \mw 1 \mw 1
\end{align*}
Thus the rank on \( \myCSW{4}{2,-2} \) is 2.  

The boundary image of \( \myCSW{5}{2,-2} \) becomes 0 as shown below. 
\begin{align*}
    \pdel( \ddx \mw x \ddx \mw x^2 \ddx \mw 1 \mw 1 )
      =& 
\Sbt{\ddx}{ x\ddx} \mw x^2 \ddx \mw 1 \mw 1 
    + x\ddx \mw \Sbt{\ddx}{ x^2 \ddx} \mw 1 \mw 1  \\& 
+ x\ddx \mw  x^2 \ddx \mw \Sbt{\ddx}{1}  \mw 1  
+ x\ddx \mw  x^2 \ddx \mw 1 \mw  \Sbt{\ddx}{1}  \\& 
- \ddx \mw (
\Sbt{x\ddx} { x^2 \ddx } \mw 1 \mw 1 + 
{ x^2 \ddx } \mw \Sbt{x\ddx}{1} \mw 1 + 
{ x^2 \ddx } \mw 1 \mw \Sbt{x\ddx}{1}  
) \\& 
+ \ddx \mw x \ddx \mw ( 0 \mw 1 + 1 \mw 0) + 0 \\ = & 
\ddx \mw x^2 \ddx \mw 1 \mw 1 
    + x\ddx \mw  2 { x \ddx} \mw 1 \mw 1  \\& 
+ x\ddx \mw  x^2 \ddx \mw 0   \mw 1  
+ x\ddx \mw  x^2 \ddx \mw 1 \mw  0  \\& 
- \ddx \mw (
{ x^2 \ddx } \mw 1 \mw 1 + 
{ x^2 \ddx } \mw 0  \mw 1 + 
{ x^2 \ddx } \mw 1 \mw  0   
)  
    \\ 
 = & \ddx \mw x^2 \ddx \mw 1 \mw 1 
    + 2 x\ddx \mw  { x \ddx} \mw 1 \mw 1  \\& 
+ 0 + 0   - \ddx \mw ( { x^2 \ddx } \mw 1 \mw 1 + 0 +  0   ) \\ 
 =& 
\ddx \mw x^2 \ddx \mw 1 \mw 1 
    + 2 x\ddx \mw   { x \ddx} \mw 1 \mw 1  
- \ddx \mw ( { x^2 \ddx } \mw 1 \mw 1  )  \\ 
=& 2 x\ddx \mw   { x \ddx} \mw 1 \mw 1  
= 0 \;, 
\quad \text{ and } \\
    \pdel  \myCSW{5}{2,-2}  = & 0 \;. 
\end{align*}
\kmcomment{ 
\begin{align*}
    \pdel  \myCSW{5}{2,-2}  = & 
    \pdel ( \kmMat{0}{0} \mw \kmMat{0}{1} \mw \kmMat{0}{2} ) 
 \mw  \kmMat{1}{0} \mw \kmMat{1}{0} 
 \\ =& 
 ( \kmMat{0}{0} \mw \kmMat{0}{2}  
+ 2 \kmMat{0}{1} \mw \kmMat{0}{1}  
-   \kmMat{0}{0} \mw \kmMat{0}{2} )  
 \mw  \kmMat{1}{0} \mw \kmMat{1}{0} 
= 0 
\end{align*}
}
Betti numbers are given as follows. \(
\begin{array}{c|*{4}{c}}
    & C_2 & C_3 & C_4 & C_5 \\\hline
    \text{dim} & 2 & 3 & 2 & 1 \\\hline
    \text{rank} & 0 & 1 & 2 & 0 \\\hline
    \text{Betti} & 1 & 0 & 0 & 1
\end{array}
\)

\subsection{\(n=1, w=-h=3\)} 
By Theorem \ref{thm:n:1}, \(m\) of the chain space \( \myCSW{m}{3,-3} \)  
runs \( 3,4,5,6\). 
\kmcomment{
When \(m=3\), \(i_1+i_2+i_3=w=3\) and \(k_1 + k_2 +k_3 =  h+3 = 0\) hold and
so \(k_1 = k_2 =k_3 = 0\),  and those imply 
}
\begin{align*} \myCSW{3}{3,-3}&=\kmMat{0}{0} \mw \kmMat{1}{0} \mw \kmMat{2}{0} 
+ \kmMat{1}{0} \mw \kmMat{1}{0} \mw \kmMat{1}{0}\ , \\
\myCSW{4}{3,-3}&=\kmMat{0}{0}\mw\kmMat{0}{1}\mw\kmMat{1}{0}\mw\kmMat{2}{0} 
+ \kmMat{0}{0}\mw\kmMat{1}{0}\mw\kmMat{1}{0}\mw\kmMat{1}{1} 
             \\& 
+ \kmMat{0}{1}\mw\kmMat{1}{0}\mw\kmMat{1}{0}\mw\kmMat{1}{0}\ , \\
\myCSW{5}{3,-3}&=
\kmMat{0}{0} \mw  \kmMat{0}{1} \mw   \kmMat{1}{0} \mw \kmMat{1}{0}  \mw \kmMat{1}{1} 
+ \kmMat{0}{0} \mw  \kmMat{0}{2}  \mw   \kmMat{1}{0} \mw
 \kmMat{1}{0} \mw  \kmMat{1}{0}   \ , \\
\myCSW{6}{3,-3}&=\kmMat{0}{0} \mw  \kmMat{0}{1} \mw   \kmMat{0}{2}  \mw 
\kmMat{1}{0} \mw \kmMat{1}{0} \mw  \kmMat{1}{0}  \ . 
\end{align*}
We see that  
\(    \myCSW{3}{3,-3}  = \myCSW{2}{2,-2} \mw {1}\) , 
\(    \myCSW{4}{3,-3}  = \myCSW{3}{2,-2} \mw {1}\) , 
\(    \myCSW{5}{3,-3}  = \myCSW{4}{2,-2} \mw {1}\) , and 
\(    \myCSW{6}{3,-3}   = \myCSW{5}{2,-2} \mw {1}\). 
 
\kmcomment{ 
\begin{alignat*}{2}
    \myCSW{3}{3,-3} & = \myCSW{2}{2,-2} \mw {1} , \qquad & 
    \myCSW{4}{3,-3} & = \myCSW{3}{2,-2} \mw {1} , \\ 
    \myCSW{5}{3,-3} & = \myCSW{4}{2,-2} \mw {1} , &   
    \myCSW{6}{3,-3} &  = \myCSW{5}{2,-2} \mw {1} . 
 \end{alignat*}
}
We expect that 
\( \myCSW{m+\ell}{\ell+2,-\ell-2}=\myCSW{m}{2,-2} ( \mw {1}) ^{\ell}\) 
holds for \(m=2 \ldots  5\) and \( \ell > 0 \) and 
the table of boundary
operation and so the Betti numbers should be invariant.  

\subsection{General case  under \(n=1, w=-h\)} 
A given primary weight \(w\) of 
  \(\myCSW{m}{w,-w}\) with \(n=1\), it is already known \(m \) runs
  in \( M_{w} := [ w, w+1, w+2, w+3] \) in general and   
   \( M_{0} := [  1, 2, 3] \).  
We have studied 
  \(\myCSW{m}{w,-w}\) for \( w\) from 0 to 3 and have some expectation as
  wrote above.  Here, we try to find concrete expression 
  of \(\myCSW{m}{w,-w}\) for general \(w\) above. 

Since we have handled cases of \( w \leq 3\), we now assume  \( w > 3\) and
recall a typical picture and related definitions.
\begin{align}
    & \hspace{-20mm} 
 \overbrace{\kmMat{0}{i_1} \mw \cdots \mw \kmMat{0}{i_a} }^{a} 
 \mw \overbrace{ \kmMat{1}{j_1} \mw\cdots\mw \kmMat{1}{j_b}}^{b} 
 \mw \overbrace{ \kmMat{2}{k_1} \mw \cdots\mw \kmMat{2}{k_c}}^{c}  
 \label{picDesu} \\\shortintertext{where 
 \( \{i_s \} \) and \( \{k_s \} \) are sets of distinct elements.} 
    m  &= a+b+c \, , \label{wa} \\ w &=  0 \cdot  a + 1\cdot b + 2 \cdot c =
    b+ 2c \, ,  \label{pw}
    \\
m-w &=  m +h = \sum_{s=1}^{a} i_s + \sum_{s=1}^{b} j_s + \sum_{s=1}^{c} k_s
\, . \label{sw}
 \end{align}   

 \kmcomment{
Let \(L\) be a function of 0,1 defined by 
\begin{equation}   L(x)  := \begin{cases} 0 & \text{if } x=0 \\ 
\emptyset &\text{if } x=1 \end{cases}
\qquad 
\text{ 
and  the cardinality }  
  \# L(x) = 1 - x \, . 
\end{equation}
}
Here after, we assume \(\alpha\) and/or \(\beta\) run in (0,null). If
\(\#\alpha = 0\), then
\(\alpha = null\).  If \(\#\alpha= 1\), then \(\alpha = 0\), where \(\#\alpha =
\#\{\alpha\} = \) the cardinality of the set \(\{\alpha\}\) in precise.  

\subsubsection{$m=w$}
 Let \( m = w \).  \eqref{sw} implies \( i_s = j_s = k_s =0\) and we get 4
 cases as follows. 
\[ C[1]= \overbrace{ \emptyset  }^{a=0} 
 \mw \overbrace{ \kmMat{1}{0} \mw\cdots\mw \kmMat{1}{0}}^{b} 
 \mw \overbrace{ \emptyset }^{c=0} \qquad  
 C[2] = 
 \overbrace{ \emptyset  }^{a=0} 
 \mw \overbrace{ \kmMat{1}{0} \mw\cdots\mw \kmMat{1}{0}}^{b} 
 \mw \overbrace{ \kmMat{2}{0} }^{c=1} \] 
\[C[3] = \overbrace{\kmMat{0}{0} }^{a=1} 
 \mw \overbrace{ \kmMat{1}{0} \mw\cdots\mw \kmMat{1}{0}}^{b} 
 \mw \overbrace{ \emptyset }^{c=0} \qquad  
 C[4] = 
 \overbrace{\kmMat{0}{0}  }^{a=1} 
 \mw \overbrace{ \kmMat{1}{0} \mw\cdots\mw \kmMat{1}{0}}^{b} 
 \mw \overbrace{ \kmMat{2}{0} }^{c=1} \] 
Calculating the primary weight from these pictures \(C[i]\) for \(i=1,2,3,
4\), we get \(m, m-1, m+1, m\). Thus, we conclude the chain space is given
by \( C[1] + C[4]\).  
Using \(\alpha\) and \(\beta\),  we get a common table as follows: 
\[ \begin{array}{c|c|c|c}
\text{\almostprimary} & \text{0-series} &\text{ 1-series} & \text{2-series } \\ \hline
\text{almost 2ndary} & \alpha & 0,\ldots,0 & \beta \\ \hline
\text{cardinality} & \#\alpha & m - \#\alpha - \#\beta & \#\beta
\end{array} \] 
This tells that the primary weight is \( 0\cdot 
 \#\alpha + 1 \cdot ( m - \#\alpha - \#\beta) + 2 \cdot   \#\beta
 = m - \#\alpha + \#\beta \). 
 Since the primary weight should be \(m\), we get  
 \(\#\alpha = \#\beta \) and  \( C[1] + C[4]\).  
Namely, 
\begin{align}
    \myCSW{w}{w,-w} & = C[1] + C[4]=  
  \overbrace{ \kmMat{1}{0} \mw\cdots\mw \kmMat{1}{0}}^{w} 
 + 
{\kmMat{0}{0}  }
 \mw \overbrace{ \kmMat{1}{0} \mw\cdots\mw \kmMat{1}{0}}^{w-2} 
 \mw  \kmMat{2}{0} 
\\  
&=  \overbrace{ \kmMatt{1}{0} \mw\cdots\mw \kmMatt{1}{0}}^{w-2} 
\mw \left( 
 \kmMat{1}{0} \mw \kmMat{1}{0} + 
 \kmMat{0}{0} \mw  \kmMat{2}{0} \right)
\\  
&=  \overbrace{ \kmMatt{1}{0} \mw\cdots\mw \kmMatt{1}{0}}^{w-2} 
    \mw \myCSW{2}{2,-2} \, .
\end{align}

\subsubsection{\(m=w+1\)}
Let \( m= w+1\). The secondary weight condition  \eqref{sw} says that one
sum is 1 and the other sums are 0.  

Assume first that  
\( \sum_{s=1}^{a} i_s =1 , \sum_{s=1}^{b} j_s = \sum_{s=1}^{c} k_s = 0\). 
Then in \eqref{picDesu},  a-zone has two possibilities of \( \kmMat{0}{1} \) or 
    \( \kmMat{0}{1} \mw \kmMat{0}{0} \) and c-zone has possibly two cases  
  \( \kmMat{2}{0} \) or \(\emptyset\), and we get 4
 possibilities as one table.  
\[ \begin{array}{c|c|c|c}
\text{\almostprimary} & \text{0-series} &\text{ 1-series} & \text{2-series } \\ \hline
\text{almost 2ndary} & 1, \alpha & 0,\ldots,0 & \beta \\ \hline
\text{cardinality} & 1+ \#\alpha & m -1- \#\alpha - \#\beta & \#\beta
\end{array} \] 
The primary weight condition tells that 
 \(\#\alpha = \#\beta \) 
and  we get 
\begin{align*}
{\kmMat{0}{1}}
 \mw \overbrace{ \kmMat{1}{0} \mw\cdots\mw \kmMat{1}{0}}^{w} 
 + 
{\kmMat{0}{0}} \mw {\kmMat{0}{1}} \mw 
\overbrace{ \kmMat{1}{0} \mw\cdots\mw \kmMat{1}{0}}^{w-2} \mw \kmMat{2}{0}
\end{align*}
Next assume that  
\( \sum_{s=1}^{a} i_s =0 , \sum_{s=1}^{b} j_s =1,  \sum_{s=1}^{c} k_s = 0\). 
Then in \eqref{picDesu},  a-zone has two possibilities of \( \kmMat{0}{0},
\emptyset \)  and 
 c-zone has two  
 possibilities of \( \kmMat{2}{0},  \emptyset\), and we get 4
 possibilities.  We evaluate the primary weight for those possibilities, 
\[ \begin{array}{c|c|c|c}
\text{\almostprimary} & \text{0-series} &\text{ 1-series} & \text{2-series } \\ \hline
\text{almost 2ndary} &  \alpha & 1,0,\ldots,0 & \beta \\ \hline
\text{cardinality} &  \#\alpha & m - \#\alpha - \#\beta & \#\beta
\end{array} \] 
The table above and the primary weight condition tells that \( m-1 = m 
  - \#\alpha - \#\beta + 2 \#\beta \), i.e.,  
  \( \#\beta = 1+ \#\beta \), i.e.,  \(\alpha=0\) and \(\beta =
  \emptyset\), so we get 
\begin{align*}
\kmMat{0}{0} 
 \mw \overbrace{ \kmMatt{1}{0} \mw\cdots\mw \kmMatt{1}{0}}^{w-1} \mw  
\kmMat{1}{1} 
\end{align*}
Lastly, assume  that  
\( \sum_{s=1}^{a} i_s = \sum_{s=1}^{b} j_s =0,  \sum_{s=1}^{c} k_s = 1\). 
Then in \eqref{picDesu},  a-zone has two possibilities of \( \kmMat{0}{0},
\emptyset \)  and 
 c-zone has two  
 possibilities of \( \kmMat{2}{1}\) or 
 \( \kmMat{2}{1} \mw \kmMat{2}{0}\), 
 possibilities.  
In short, we use the magic? table as follows.  
\[ \begin{array}{c|c|c|c}
\text{\almostprimary} & \text{0-series} &\text{ 1-series} & \text{2-series } \\ \hline
\text{almost 2ndary} &  \alpha & 0,\ldots,0 & 1, \beta \\ \hline
\text{cardinality} &  \#\alpha & m - 1- \#\alpha - \#\beta & 1+\#\beta
\end{array} \] 
This time,  the primary weight condition shows that 
\(\#\alpha - \#\beta =2 \)  and there is no solution.

Combining those, we conclude that 
\begin{align}
\myCSW{w+1}{w,-w} & =
\kmMat{0}{1}  \mw \overbrace{ \kmMat{1}{0} \mw\cdots\mw \kmMat{1}{0}}^{w} 
+ \kmMat{0}{0}  \mw  \kmMat{0}{1}\mw 
\overbrace{ \kmMat{1}{0} \mw\cdots\mw \kmMat{1}{0}}^{w-2} 
 \mw  \kmMat{2}{0} \\& 
 \quad + \kmMat{0}{0} \mw 
   \overbrace{ \kmMat{1}{0} \mw\cdots\mw \kmMat{1}{0}}^{w-1} \mw  \kmMat{1}{1}  \notag
 \\ &= 
 \overbrace{ \kmMat{1}{0} \mw\cdots\mw \kmMat{1}{0}}^{w-2} \mw 
 ( \kmMat{0}{1}  \mw  \kmMat{1}{0}  \mw  \kmMat{1}{0} 
 + \kmMat{0}{0}  \mw  \kmMat{0}{1}  \mw  \kmMat{2}{0} 
 \\ & \qquad 
 + \kmMat{0}{0} \mw   \kmMat{1}{0} \mw   \kmMat{1}{1}  ) 
 \notag
  \\&= 
 \overbrace{ \kmMatt{1}{0} \mw\cdots\mw \kmMatt{1}{0}}^{w-2} \mw 
    \myCSW{3}{2,-2} \, . 
\end{align}

\kmcomment{
Let \(L\) be a function of 0,1 defined by 
\begin{equation}   L(x)  := \begin{cases} 0 & \text{if } x=0 \\ 
\emptyset &\text{if } x=1 \end{cases}
\qquad 
\text{ 
and  the carnality }  
  \# L(x) = 1 - x \, . 
\end{equation}
}

\subsubsection{\(m=w+2\)} \label{subsub:m:2}
Applying \(m-w=2\) to the three sets 
\(  \{i_s \}_{s=1}^{a},   \{j_s \}_{s=1}^{b},   \{k_s \}_{s=1}^{c} \)
where \(m-w = \sum_{s=1}^{a} i_s + \sum_{s=1}^{b} j_s + \sum_{s=1}^{c} k_s \), 
0 is distributed to one set and 1 to others, or 2 for one set  and 0 to others.  Thus, there are
six cases 
(2,0,0), (0,2,0),  (0,0,2), 
(1,1,0), (1,0,1),  (0,1,1) 
for  \( (\sum_{s=1}^{a} i_s , \sum_{s=1}^{b} j_s , \sum_{s=1}^{c} k_s )\).  

For the first trial, take (2,0,0), i.e., 
\( \sum_{s=1}^{a} i_s =2, \sum_{s=1}^{b} j_s =0 , \sum_{s=1}^{c} k_s=0 \).
This implies that 0-series is (2) or (0,2), and \(j_s=0\) and \(k_s=0\). 
We remember that 0-series and 2-series are skew-symmetric and (1,1, \ldots) is not allowed.    
Now the picture  \eqref{picDesu} becomes  in short 
\[ \begin{array}{c|c|c|c}
\text{\almostprimary} & \text{0-series} &\text{ 1-series} & \text{2-series } \\ \hline
\text{almost 2ndary} & 2, \alpha & 0,\ldots,0 & \beta \\ \hline
\text{carnality} & 1+ \#\alpha & m - 1- \#\alpha - \#\beta & \#\beta
\end{array} \] 
This time, the primary weight condition is 
\( m-2 = w = 1\cdot(m-1 -\#\alpha  - \# \beta ) + 2 \# \beta = 
m -1 - \#\alpha + \# \beta \) and we get 
\( \#\alpha = 1 + \# \beta \), i.e.,  
\( \#\alpha  = 1 ,  \# \beta = 0 \) more exactly  
\( \alpha  = 0 ,   \beta = \emptyset  \). 

For the second trial,  take (1,1,0), i.e., 
\[ \begin{array}{c|c|c|c}
\text{\almostprimary} & \text{0-series} &\text{ 1-series} & \text{2-series } \\ \hline
\text{almost 2ndary} & 1, \alpha & 1,0,\ldots,0 & \beta \\ \hline
\text{carnality} & 1+ \#\alpha & m - 1- \#\alpha - \#\beta & \#\beta
\end{array} \] 
Checking the primary weight condition again, we get 
\( \alpha  = 0 ,   \beta = \emptyset  \). 

We continued the completely same the discussion and find no more
contribution to the chain space.  
We do not show those discussion here but leave the detail in Appendix. 

Our conclusion here is as follows. 
\begin{align}
    \myCSW{w+2}{w,-w} & =
\kmMat{0}{0}  \mw 
\kmMat{0}{2}  \mw 
\overbrace{ \kmMat{1}{0} \mw\cdots\mw \kmMat{1}{0}}^{w} 
+ \kmMat{0}{0}  \mw  \kmMat{0}{1}\mw 
\overbrace{ \kmMat{1}{0} \mw\cdots\mw \kmMat{1}{0}}^{w-1} 
 \mw  \kmMat{1}{1} 
 \\ &= 
 \overbrace{ \kmMat{1}{0} \mw\cdots\mw \kmMat{1}{0}}^{w-1} \mw 
 ( \kmMat{0}{0}  \mw  \kmMat{0}{2}  \mw  \kmMat{1}{0} 
 + \kmMat{0}{0}  \mw  \kmMat{0}{1}  \mw  \kmMat{1}{1} ) 
  \\&= 
 \overbrace{ \kmMatt{1}{0} \mw\cdots\mw \kmMatt{1}{0}}^{w-2} \mw 
    \myCSW{4}{2,-2} \, . 
\end{align}
\subsubsection{\(m=w+3\)} \label{w=m:3}
Applying \(m-w=3\) to the three sets 
\( \{ i_s \}_{s=1}^{a},  \{ j_s \}_{s=1}^{b}, \{ k_s\}_{s=1}^{c} \) 
where \(m-w = \sum_{s=1}^{a} i_s + \sum_{s=1}^{b} j_s + \sum_{s=1}^{c} k_s \), 
and 3 cases: the first one is  
\(\sum_{s=1}^{a} i_s = \sum_{s=1}^{b} j_s = \sum_{s=1}^{c} k_s =1 \),  
the second one is  
 \(\{\sum_{s=1}^{a} i_s , \sum_{s=1}^{b} j_s , \sum_{s=1}^{c} k_s =1 \}
={0,1,2} \} \) as a set, and the last case is that one is distributed by 3 and
the others are distributed by 0.   
We introduce here the case which give us the solution, the rest will be
reported in Appendix, too. 

Suppose that 0-series has the number 3. Then we have a table (joined two). 
\kmcomment{
\begin{align*} 
  &  \begin{array}{c|c|c|c}
\text{\almostprimary} & \text{0-series} &\text{ 1-series} & \text{2-series } \\ \hline
\text{almost 2ndary} & 3, \alpha & 0,\ldots,0 & \beta \\ \hline
\text{carnality} & 1+ \#\alpha & m - 1- \#\alpha - \#\beta & \#\beta
\end{array} 
\\
  &   \begin{array}{c|c|c|c}
\text{\almostprimary} & \text{0-series} &\text{ 1-series} & \text{2-series } \\ \hline
\text{almost 2ndary} & 1,2, \alpha & 0,\ldots,0 & \beta \\ \hline
\text{carnality} & 2+ \#\alpha & m - 2- \#\alpha - \#\beta & \#\beta
\end{array} 
\end{align*} 
}

\begin{align*} 
    \begin{array}{c|c|c|c|*{3}{|c}}
        \hline
\text{\almostprimary} & \text{0-series} &\text{ 1-series} & \text{2-series }
                      & 
 \text{0-series} &\text{ 1-series} & \text{2-series } \\ \hline 
\text{almost 2ndary} & 3, \alpha & 0,\ldots,0 & \beta
                     & 
 1,2, \alpha & 0,\ldots,0 & \beta \\ \hline
\text{carnality} & 1+ \#\alpha & m - 1- \#\alpha - \#\beta & \#\beta
 & 2+ \#\alpha & m - 2- \#\alpha - \#\beta & \#\beta \\\hline
\end{array} 
\end{align*}
We check the primary weight condition.  

Case 1. \( m-3 = w = m -1 - \# \alpha + \# \beta\) and 
\( \alpha  = 2 +  \beta  \) and it is impossible.

Case 2. \( m-3 = w = m -2 - \# \alpha + \# \beta \) and 
\( \alpha  = 1 +  \beta\), namely  \(\alpha = 0\) and \(\beta=\emptyset\).  

About the rest of candidates, we  will show that there is no solution in Appendix
and conclude as follows.  
\begin{align}
    \myCSW{w+3}{w,-w}  =&
\kmMat{0}{0}  \mw 
\kmMat{0}{1}  \mw 
\kmMat{0}{2}  \mw 
\overbrace{ \kmMat{1}{0} \mw\cdots\mw \kmMat{1}{0}}^{w} 
= \overbrace{ \kmMatt{1}{0} \mw\cdots\mw \kmMatt{1}{0}}^{w-2} \mw 
 \myCSW{5}{2,-2} \, . 
\end{align}
Combining the above results, we have an equation. 
\begin{theorem} \label{thm:II}
    For \( w > 2\), we have the next relations. 
\begin{equation}
\left[
    \begin{array}{c} 
    \myCSW{w+0}{w,-w} \\ 
    \myCSW{w+1}{w,-w} \\
    \myCSW{w+2}{w,-w} \\
    \myCSW{w+3}{w,-w} \end{array}
\right] =  
 \overbrace{ \kmMatt{1}{0} \mw\cdots\mw \kmMatt{1}{0}}^{w-2} \mw 
 \left[
    \begin{array}{c} 
    \myCSW{2+0}{2,-2} \\ 
    \myCSW{2+1}{2,-2} \\
    \myCSW{2+2}{2,-2} \\
    \myCSW{2+3}{2,-2} \end{array}
\right]
    \end{equation}
\end{theorem}
By Theorem \ref{thm:II} above, we expect commutativity of the following diagram. 
\begin{align}
    \begin{CD}
\myCSW{w+0}{w,-w} @<\pdel<<  
    \myCSW{w+1}{w,-w}  @<\pdel <<  
    \myCSW{w+2}{w,-w}  @<\pdel <<  
    \myCSW{w+3}{w,-w} 
    \\
    @A{\I}AA @A{\I}AA @A{\I}AA @A{\I}AA  \\
    \myCSW{2+0}{2,-2} @<\pdel << 
    \myCSW{2+1}{2,-2} @<\pdel <<
    \myCSW{2+2}{2,-2} @<\pdel <<
    \myCSW{2+3}{2,-2}
    \\
        \end{CD}
        \label{CD:one}
\end{align} where \(\I\) is the multiplication operator by 
 \(\overbrace{ 1 \mw\cdots\mw 1 }^{w-2} \in  
 \overbrace{ \kmMat{1}{0} \mw\cdots\mw \kmMat{1}{0}}^{w-2}  \). 
 \kmcomment{  
 We denote general target by 
 \(\overbrace{ 1 \mw\cdots\mw 1 }^{p} 
 \mw A_1 \mw \cdots \mw A_a \mw   
 \mw B_1 \mw \cdots \mw B_b \mw   
 \mw C_1 \mw \cdots \mw C_c \) and 
}

\( \Sbt{1}{ \kmMat{0}{j}} = 0 \) and 
\( \Sbt{1}{ \kmMat{2}{j}} = 0 \) for every  j. 
Since 
\( \kmMat{1}{j} = \mR x^j  \),  
\( \Sbt{1}{ \kmMat{1}{j}} = \kmMat{2}{j-1} \), and 
\( \Sbt{1}{ \kmMat{1}{0}} = 0  \) in particular. 
\kmcomment{ 
 \( \Sbt{1}{ A_j } = d( A_j )\) for differential forms, and
it vanishes for 1-form because of \(n=1\). 
For 0-form \(f\), 
\( \Sbt{1}{ f } = d( f ) = \frac{ d f}{ d x} d x \).  
} 

To understand relations of \( \pdel (\I \mw\text{object})\) and \( \I
\mw ( \pdel ( \text{object})) \), we have to deal with element of the next
form. 
\begin{equation}
 \overbrace{\kmMatt{1}{0} \mw \cdots \mw \kmMatt{1}{0} }^{w-2} 
 \mw \overbrace{\kmMat{0}{i_1} \mw \cdots \mw \kmMat{0}{i_a} }^{a} 
 \mw \overbrace{ \kmMat{1}{j_1} \mw\cdots\mw \kmMat{1}{j_b}}^{b} 
 \mw \overbrace{ \kmMat{2}{k_1} \mw \cdots\mw \kmMat{2}{k_c}}^{c}  
 \label{picDesu:2} \end{equation} where 
 \( \{i_s \} \) and \( \{k_s \} \) are sets of distinct elements. 

A general element of the space above is given by  
 \(\overbrace{ 1 \mw\cdots\mw 1 }^{w-2} 
 \mw A_1 \mw \cdots \mw A_a    
 \mw B_1 \mw \cdots \mw B_b    
 \mw C_1 \mw \cdots \mw C_c \) where 
 \( A_s \in \kmMat{0}{i_s},  
 B_s \in \kmMat{1}{j_s},   C_s \in \kmMat{2}{k_s}\).   
 We often denote 
 \(A_1 \mw \cdots \mw A_a \) by \(\AAs\), and \(\BBs\) and \(\CCs\) also. 
We easily split 
   \( \pdel ( \I \mw \AAs \mw \BBs \mw \CCs ) \) to two parts. 
\begin{equation}
    \pdel ( \I \mw \AAs \mw \BBs \mw \CCs )  = \text{no}{\I}{\text
    contrib} + {\I}\text{contrib} = \parity{w-2} \I \mw \pdel (\AAs 
    \mw \BBs  \mw \CCs ) + {\I}\text{contrib} \end{equation}
\begin{align}
{\I}\text{contrib}  & = \partt[0] + \partt[1] + \partt[2] +
    \partt[3], \text{ where }  \\
\partt[0] & := \sum_{i<j\leq w-2} \parity{i-1 +  \sum_{i<s<j} 1 }
1 \mw  \cdots \mw \Sbt{1}{1} \mw \cdots \mw 1 \mw \AAs \mw \BBs \mw \CCs = 0 
\\ 
    \partt[1] &:=  \sum_{i \leq w-2, j \leq a } \parity{ i-1 + (w-2-i) + }
    1^{w-2 - 1} \mw A_{i}\mw \cdots \mw \Sbt{1}{ A_{i}} \mw \cdots \mw
    \BBs \mw \CCs = 0 \\
    \partt[2] &:= \sum_{i\leq w-2, j \leq b} \parity{i-1 + w-2 -i +j-1 
    } 1^{w-3} \mw \AAs \mw \cdots \mw \Sbt{1}{B_j}\mw \cdots \mw B_b \mw  \CCs  \\ 
 &= \parity{  w- 2 + b }\sum_{ j \leq b}  1^{w-3} \mw \AAs \mw 
 \cdots \widehat{B_j} \cdots \mw B_b \mw \Sbt{1}{B_{j}} \mw  \CCs \; \notag \\ 
 &= 
 \parity{  w- 2  }(1\mw)^{w-3} \mw \AAs \mw 
 ( \ActOne  (\BBs) ) 
 \mw  \CCs \;,  \notag 
\\\shortintertext{where} 
        \ActOne(\BBs) & := \begin{cases}  \parity{b} \sum_{j \leq b} B_1 \mw 
 \cdots \widehat{B_j} \cdots \mw B_b \mw \Sbt{1}{B_{j}} \quad \text{for } b>
 0 \; \\ 0 \quad   \text{ otherwise. }  \end{cases}
\\
\partt[3] &:= \sum_{k\leqq c} \parity{ \text{omit detail} } 1^{w-3} \mw \AAs
\mw \BBs 
\mw C_1 \mw \cdots \mw \Sbt{1}{C_k} \mw \cdots = 0 \; . \\
\kmcomment{Full Expression 
    \partt[3] &:= \sum_{k\leqq c} \parity{i-1 + (w-2-i) + 0 + b + 2\cdot
    (k-1)  } 1^{w-3} \mw \AAs \mw \BBs \mw C_1 \mw \cdots \mw \Sbt{1}{C_k} \mw
    \cdots = 0\; .  \notag\\
}
\shortintertext{Thus, we get
    \( \I\text{contrib} = \partt[2]\), namely} 
        \I\text{contrib} &= 
        \parity{  w- 2  } 1^{w-3} \mw \AAs \mw  \ActOne  (\BBs) \mw \CCs   
 \end{align}
Finally, we combine the several formulae above and get a result.   
\begin{theorem} \label{thm:dAndII}
Under the same notations above, we have 
\begin{equation}
       ( \pdel \circ ( \I \mw ) )  (\AAs \mw \BBs  \mw \CCs )  = 
\parity{  w- 2  }(1\mw)^{w-3} \mw \Big( 
 1  \mw  \pdel (\AAs \mw \BBs  \mw \CCs )  +   
 \AAs \mw  \ActOne ( \BBs)  \mw  \CCs \Big)  \;  
 \label{thm5:part0}
\end{equation}
for \(  \AAs \mw \BBs  \mw \CCs \in 
    \myCSW{2+k}{2,-2} \) for \(k=0,1,2,3\) and for \( w >2 \). 
\end{theorem}
\begin{remark}
Although 
we expected commutativity of the diagram \eqref{CD:one}, it does not hold in
general.  
However, we get some kind of shortcut as follows. 
\begin{equation}
       \left( \pdel \circ ( \I \mw )   - 
\parity{w-2}( \I \mw )\circ \pdel \right) (\AAs \mw \BBs  \mw \CCs )  =   
\parity{  w- 2  }(1\mw)^{w-3} \mw \AAs \mw  \ActOne  ( \BBs)  \mw  \CCs \;  
 \label{thm5:part1}
\end{equation}
for \(  \AAs \mw \BBs  \mw \CCs \in 
    \myCSW{2+k}{2,-2} \) for \(k=0,1,2,3\) and for \( w >2 \). 
\end{remark}



\kmcomment{ 
Does II (the operator we defined) work like as the Euler vector field? 

Now we have 
\( \pdel ( \I \mw A_{m+1} \mw \cdots \mw A_{m+r} ) 
= \parity{m} \I \mw \pdel ( \mw A_{m+1} \mw \cdots \mw A_{m+r} ) +\partt[2] \)

We NEED some lucky punch.  

We study about \(\partt[2]\) in detail using  the picture 
 \(\overbrace{ 1 \mw\cdots\mw 1 }^{p} \mw   
\overbrace{\texttt{  1-vecs} }^{a} \mw
\overbrace{ \texttt{ 0-forms} }^{b} \mw
\overbrace{\texttt{ 1-forms} }^{c} 
 \). We denote general target by 
 \(\overbrace{ 1 \mw\cdots\mw 1 }^{p} 
 \mw A_1 \mw \cdots \mw A_a 
 \mw B_1 \mw \cdots \mw B_b 
 \mw C_1 \mw \cdots \mw C_c \) and 
\kmcomment{Old style: 
    \partt[2] &:= \sum_{i \leq p<j} \parity{ i-1 + (p-i) + \sum_{p<s<j} a_{s}}
    (1\mw)^{p - 1} \mw A_{p+1}\mw \cdots \mw \Sbt{1}{ A_{j}} \mw \cdots \\
}
\begin{align*}
    \partt[2] &:= 0 (\text{for 1-vecs}) \\
              &  + \sum_{i \leq p,  j \leq b} \parity{ i-1 + (p-i) + j-1 }
              ( 1\mw)^{p - 1} \mw A_{1}\mw \cdots \mw A_a \mw B_1 \mw \cdots \Sbt{1}{ B_{j}}
    \mw \cdots \mw B_{b} \mw \CCs 
            \\& 
+  0 (\text{for 1-forms})  
    \\ &=  \sum_{i \leq p,  j \leq b} \parity{ i-1 + (p-i) + j-1+b-j }
              (1\mw)^{p - 1} \mw A_{1}\mw \cdots \mw A_a \mw B_1 \mw \cdots 
    \widehat{B_j}  \cdots \mw B_{b} \mw \Sbt{1}{ B_{j}}\mw \CCs 
    \\ &= \parity{ p + b } \sum_{i \leq p,  j \leq b} 
    (1\mw)^{p - 1} \mw A_{1}\mw \cdots \mw A_a \mw B_1 \mw \cdots 
    \widehat{B_j} \cdots \mw B_{b} \mw \Sbt{1}{ B_{j}}\mw \CCs 
    \\ &= p     \parity{ p + b }
    ( 1 \mw) ^{p - 1} \mw A_{1}\mw \cdots \mw A_a \mw 
    \sum_{j\leq b} 
    B_1 \mw \cdots 
    \widehat{B_j} \cdots \mw B_{b} \mw \Sbt{1}{ B_{j}}\mw \CCs 
\end{align*}
}

\kmcomment{
\begin{align}
    \myCSW{w}{w,-w} & = C[1] + C[4]=  
  \overbrace{ \kmMat{1}{0} \mw\cdots\mw \kmMat{1}{0}}^{w} 
 + 
{\kmMat{0}{0}  }
 \mw \overbrace{ \kmMat{1}{0} \mw\cdots\mw \kmMat{1}{0}}^{w-2} 
 \mw  \kmMat{2}{0} 
}
  

\kmcomment{
 \overbrace{ \emptyset  }^{a} 
 \mw \overbrace{ \kmMat{1}{0} \mw\cdots\mw \kmMat{1}{0}}^{b} 
 \mw \overbrace{ \emptyset }^{c} 
 + 
 \overbrace{\kmMat{0}{0}  }^{a} 
 \mw \overbrace{ \kmMat{1}{0} \mw\cdots\mw \kmMat{1}{0}}^{b} 
 \mw \overbrace{ \kmMat{2}{0} }^{c} 
}

\kmcomment{
\[ \overbrace{\kmMat{0}{i_1} \mw \cdots \mw \kmMat{0}{i_a} }^{a} 
 \mw \overbrace{ \kmMat{1}{j_1} \mw\cdots\mw \kmMat{1}{j_b}}^{b} 
 \mw \overbrace{ \kmMat{2}{k_1} \mw \cdots\mw \kmMat{2}{k_c}}^{c} \] 

\[ \overbrace{\kmMat{0}{i_1} \mw \cdots \mw \kmMat{0}{i_a} }^{a} 
 \mw \overbrace{ \kmMat{1}{j_1} \mw\cdots\mw \kmMat{1}{j_b}}^{b} 
 \mw \overbrace{ \kmMat{2}{k_1} \mw \cdots\mw \kmMat{2}{k_c}}^{c} \] 
}

\bigskip 

\kmcomment{ 
\begin{tikzcd}
    \myCSW{w+0}{w,-w}  \ar[r,"f"] &   
    \myCSW{w+1}{w,-w}  \ar[l, "\pdel"] \\
    X \ar[u] & Y \ar[u]
\end{tikzcd}

\begin{center}
\begin{tikzcd}
    \myCSW{w+0}{w,-w} &  
    \myCSW{w+1}{w,-w}  \ar[l,"\pdel" above] &  
    \myCSW{w+2}{w,-w}  \ar[l,"\pdel" above] & 
    \myCSW{w+3}{w,-w}  \ar[l,"\pdel" above] & 
    \\
    \myCSW{2+0}{2,-2}  \ar[u,"\I" right]  & 
    \myCSW{2+1}{2,-2}  \ar[l,"\pdel"]  \ar[u,"\I" right]  &   
    \myCSW{2+2}{2,-2}  \ar[l,"\pdel"]  \ar[u,"\I" right]  & 
    \myCSW{2+3}{2,-2} \ar[l,"\pdel"]   \ar[u,"\I" right]   
    \\
\end{tikzcd}
\end{center}
}

We define new operators from the right hand side of 
 \eqref{thm5:part0} as 
\begin{equation}
 \myCSW{2+k}{2,-2} \stackrel{\pdel '}{\longleftarrow}
    \myCSW{2+k+1}{2,-2} \  : \  
        \pdel '  (\AAs \mw \BBs  \mw \CCs ) :=   
 1  \mw  \pdel (\AAs \mw \BBs  \mw \CCs )  +   
 \AAs \mw  \ActOne ( \BBs)  \mw  \CCs \; .  
\end{equation}

\begin{lemma} 
    The ranks  \(1 \mw \pdel \) and \( \pdel '\) of 
 \( \myCSW{2+k}{2,-2} \stackrel{ }{\longleftarrow}
    \myCSW{2+k+1}{2,-2} \) are  equal for \(k=0,1,2,3\). 
\end{lemma}
\textbf{Proof:}
We have to compare each rank   
of the two operators \(\pdel\) and \(\pdel '\). 
We show this basic process by tables below: 
\kmcomment{ 
\( \begin{array}{*{3}{|c}|}  
    \hline
\AAs \mw \BBs  \mw \CCs  
& \pdel (\AAs \mw \BBs  \mw \CCs )  &  \AAs \mw \ActOne (\BBs) \mw \CCs
 \\\hline  
    \text{first} & & \\\hline 
 \text{others} && \\\hline
 \end{array}
\). 
}

\medskip

\( \myCSW{5}{2,-2}\) 
 is spanned by \( \ddx \mw x \ddx \mw x^2 \ddx \mw 1 \mw 1 \) and we get 
the next table, and so on. 
\begin{align*}
& \begin{array}[c]{*{4}{|c}|}  
    \hline
    & \AAs \mw \BBs  \mw \CCs  
& \pdel (\AAs \mw \BBs  \mw \CCs )  &  \AAs \mw \ActOne (\BBs) \mw \CCs
 \\\hline \hline
\myCSW{5}{2,-2}  
&
\ddx \mw x \ddx \mw x^2 \ddx \mw 1 \mw 1 
& 0 & 0 \\\hline\hline
\myCSW{4}{2,-2}  &
\ddx \mw x \ddx \mw 1 \mw x 
& x \ddx \mw 1 \mw 1 + \ddx \mw x\ddx  \mw d x & 
  \ddx \mw x \ddx \mw 1   \mw d x
\\\hline
 & 
\ddx \mw x^2 \ddx \mw 1 \mw 1 & 2 x \ddx \mw 1 \mw 1 & 0   \\\hline\hline
\myCSW{3}{2,-2} 
   & 
\ddx \mw x \ddx \mw  d x & 0 & 0 \\\hline 
  & 
\ddx \mw 1 \mw x &  1 \mw 1 - \ddx \mw  d x  & \ddx \mw 1 \mw d x \\\hline
                 & 
\scaleto{ x \ddx \mw 1 \mw 1 }{18pt} &  0 & 0   \\\hline\hline
\myCSW{2}{2,-2}  
& \ddx \mw d x & 0 & 0 \\\hline 
& 1 \mw 1 &  0   &  0 \\\hline\hline
\end{array}
\end{align*}  

\bigskip 

Thus, we check that each rank of  \( \pdel: \myCSW{2+k}{2,-2} \leftarrow
\myCSW{2+k+1}{2,-2} \) is equal with that of \( \pdel '\).  
\qed

The Lemma above implies the next Corollary. 

\begin{Cor} For \(k=0,1,2,3\) and for \( w>2\),  it holds that 
    the rank of \( \myCSW{w+k}{w,-w} \stackrel{\pdel}{\longleftarrow}
    \myCSW{w+k+1}{w,-w} \) is equal to 
    the rank of \( \myCSW{2+k}{2,-2} \stackrel{\pdel}{\longleftarrow}
    \myCSW{2+k+1}{2,-2} \). 
\end{Cor}

%
\kmcomment{ 
\begin{align*}
    \myCSW{2}{2,-2}=& \kmMat{0}{0} \mw  \kmMat{2}{0} + 
 \kmMat{1}{0} \mw  \kmMat{1}{0}\; , \\  
    \myCSW{3}{2,-2}=& \kmMat{0}{0}\mw\kmMat{0}{1}\mw \kmMat{2}{0}
+ \kmMat{0}{0}\mw\kmMat{1}{0}\mw \kmMat{1}{1}
+ \kmMat{0}{1}\mw\kmMat{1}{0}\mw \kmMat{1}{0}
\\
    \myCSW{4}{2,-2}=&\kmMat{0}{0}\mw\kmMat{0}{1}\mw \kmMat{1}{0}\mw \kmMat{1}{1}
+ \kmMat{0}{0}\mw\kmMat{0}{2}\mw \kmMat{1}{0}\mw \kmMat{1}{0}
\\
    \myCSW{5}{2,-2}=&\kmMat{0}{0}\mw\kmMat{0}{1}\mw \kmMat{0}{2}\mw \kmMat{1}{0}
\mw  \kmMat{1}{0}\, .  \end{align*}
}

\nocite{Mik:Miz:super3}


\renewcommand{\appendixname}{Appendix } 

\kmcomment{
「appendix A」。}



{
\appendix 



\section{Full proofs for case \(m=w+2\)}
\kmcomment{ 
\titlespacing{\subsection}{0pt}{\medskipamount}{*1.5}
\titlespacing{\subsubsection}{0pt}{\medskipamount}{*0.5}
\titlespacing*{\section}
 {0pt}{5.5ex plus 1ex minus .2ex}{4.3ex plus .2ex}
 \titlespacing*{\subsection}
 {0pt}{5.5ex plus 1ex minus .2ex}{4.3ex plus .2ex}
 \titlespacing*{\subsubsection}
 {0pt}{5.5ex plus 1ex minus .2ex}{4.3ex plus .2ex}
}

Applying \(m-w=2\) to the three sets \(  \{i_s \}_{s=1}^{a}, \{j_s \}_{s=1}
^{b},  \{k_s \}_{s=1}^{c} \) where \(m-w = \sum_{s=1}^{a} i_s + \sum_{s=1}
^{b} j_s + \sum_{s=1}^{c} k_s \), 0 is distributed to one set and 1 to
others, or 2 for one set  and 0 to others.  

Thus, there are six data sets (2, 0,0), (0,2,0), (0,0,2), (1,1,0), (1,0,1), (0,1,1) for  \( (\sum_{s=1}^{a}
i_s , \sum_{s=1}^{b} j_s , \sum_{s=1}^{c} k_s )\).

\subsection{\{2,0,0\}-type}
There are 3  data (2, 0,0), (0,2,0), (0,0,2) for  \( (\sum_{s=1}^{a}
i_s , \sum_{s=1}^{b} j_s , \sum_{s=1}^{c} k_s )\).  

\kmcomment{
\titlespacing*{\chapter}
 {0pt}{0.5ex plus 1ex minus .2ex}{0.3ex plus .2ex}
\titlespacing*{\section}
 {0pt}{0.5ex plus 1ex minus .2ex}{0.3ex plus .2ex}
\titlespacing*{\subsection}
 {0pt}{0.5ex plus 1ex minus .2ex}{0.1ex plus .1ex}
\titlespacing*{\subsubsection}
{0pt}{0.5ex plus 1ex minus .2ex}{0.3ex plus .2ex}
}

\subsubsection{(2,0,0)}
\begin{wrapfigure}[3]{r}[5mm]{0.55\textwidth}
\vspace{-1.5 \baselineskip}
\( \begin{array}{c|c|c|c}
\text{\almostprimary} & \text{0-series} &\text{ 1-series}&\text{2-series }\\\hline
\text{almost 2ndary} & 2, \alpha & 0,\ldots,0 & \beta \\ \hline
\text{carnality} & 1+ \#\alpha & m - (1+ \#\alpha) - \#\beta & \#\beta
\end{array} \) 
\end{wrapfigure}
For the first trial,  take (2,0,0), i.e., 
\( \sum_{s=1}^{a} i_s =2, \sum_{s=1}^{b} j_s =0 , \sum_{s=1}^{c} k_s=0 \).
This implies that 0-series is (2,0,\ldots) or (1,1,\ldots), and \(j_s=0\) and 
\(k_s=0\). 0-series and 2-series are skew-symmetric, (1,1, \ldots) is not allowed.    
Now \eqref{picDesu:2} becomes 
The primary weight becomes \(m - (1+ \#\alpha) + \#\beta = m - 1-
\#\alpha +  \#\beta \)
 and we get \(0 = 1- \#\alpha+  \#\beta\). 
This shows that 
 \( \#\alpha = 1 + \#\beta\) and so 
 \( \#\alpha = 1,  \#\beta =0\), i.e.,  \( \alpha  =0\) and 
 \( \beta   =\emptyset\).  
 This is part of subsubsection \ref{subsub:m:2}. 

\subsubsection{(0,2,0)}
\begin{wrapfigure}[3]{r}[5mm]{0.55\textwidth}
\vspace{-1.0 \baselineskip}
\( \begin{array}{c|c|c|c}
\text{\almostprimary} & \text{0-series} &\text{ 1-series}&\text{2-series }\\\hline
\text{almost 2ndary} & 0 & 2 (\text{or } 1,1),0,\ldots,0 & \beta \\ \hline
\text{carnality} & 1 & m - 1 - \#\beta & \#\beta
\end{array} \) \end{wrapfigure}
Let \(\sum_{s=1}^{a} i_s =0, \sum_{s=1}^{b} j_s =2, \sum_{s=1}^{c} k_s=0 \).
This implies the table right.  
The primary weight becomes \(m - 1 + \#\beta\)
 and we get \(0 = 1 + \#\beta\). This says no solution.  

\subsubsection{(0,0,2)}

\begin{wrapfigure}[3]{r}[5mm]{0.55\textwidth}
\vspace{-1.5 \baselineskip}
\( \begin{array}{c|c|c|c}
\text{\almostprimary} & \text{0-series} 
                      &\text{ 1-series}&\text{2-series }\\\hline
\text{almost 2ndary} & \alpha & 0 ,0,\ldots,0 & 2, \beta \\ \hline
\text{carnality} & \#\alpha & m - \#\alpha - ( 1 + \#\beta) & 1+ \#\beta
\end{array} \) \end{wrapfigure} 
This time, we have 
The primary weight becomes \(m - \#\alpha + (1 + \#\beta ) \)
 and we get \(0 = 3 - \#\alpha + \#\beta\). This says no solution.

\subsection{\{1,1,0\}-type}
There are 3 data sets  (1,1,0), (1,0,1), (0,1,1) for  
\( (\sum_{s=1}^{a} i_s , \sum_{s=1}^{b} j_s , \sum_{s=1}^{c} k_s )\).  

\subsubsection{[1,1,0]}
\begin{wrapfigure}[3]{r}[5mm]{0.55\textwidth}
\vspace{-1.2 \baselineskip}
\( \begin{array}{c|c|c|c}
\text{\almostprimary} & \text{0-series} &\text{ 1-series}&\text{2-series }\\\hline
\text{almost 2ndary} &1, \alpha & 1 ,0,\ldots,0 & \beta \\ \hline
\text{carnality} & 1+ \#\alpha & m -(1+ \#\alpha) -  \#\beta & \#\beta
\end{array} \)  \end{wrapfigure}
Take (1,1,0) and get the  table.   
The primary weight becomes \(m -1- \#\alpha + \#\beta  \)
 and we get \(0 = 1 - \#\alpha + \#\beta\), i.e.,  \( \alpha= 0,
\beta = \emptyset\). 

\kmcomment{
We continued the completely same the discussion and find no more
contribution to the chain space.  
We do not show those discussion here but leave the detail in Appendix if
possible. }

\subsubsection{[1,0,1]}
\begin{wrapfigure}[2]{r}[5mm]{0.6\textwidth}
\vspace{-1.0 \baselineskip}
\(\begin{array}{c|c|c|c}
\text{\almostprimary} & \text{0-series} &\text{ 1-series}&\text{2-series }\\\hline
\text{almost 2ndary} &1, \alpha & 0,\ldots,0 & 1,\beta \\ \hline
\text{carnality} & 1+ \#\alpha & m -(1+ \#\alpha) - (1+ \#\beta) & 1+\#\beta
\end{array} \) 
\end{wrapfigure}
Take (1,0,1) and get the table.   
The primary weight becomes \(m - (1+ \#\alpha) + (1+ \#\beta ) \)
 and we get \(0 = 2 - \#\alpha + \#\beta\) and no solution. 

\subsubsection{[0,1,1]}
\begin{wrapfigure}[3]{r}[5mm]{0.55\textwidth}
\vspace{-1.5\baselineskip}
\(\begin{array}{c|c|c|c}
\text{\almostprimary} & \text{0-series} &\text{ 1-series}&\text{2-series }\\\hline
\text{almost 2ndary} & 0 & 1,0,\ldots,0 & 1,\beta \\ \hline
\text{carnality} & 1  & m - 1 - (1+ \#\beta) & 1+\#\beta
\end{array} \) 
\end{wrapfigure}
Take (0,1,1) and get a table.   
The primary weight becomes \(m - 1 + (1+ \#\beta ) \)
 and we get \(0 = 2 + \#\beta\) and no solution.

\section{Full proofs for case \(m=w+3\)}
Applying \(m-w=3\) to the three sets 
\( \{ i_s \}_{s=1}^{a},  \{ j_s \}_{s=1}^{b}, \{ k_s\}_{s=1}^{c} \) 
where \(m-w = \sum_{s=1}^{a} i_s + \sum_{s=1}^{b} j_s + \sum_{s=1}^{c} k_s \), 
and 3 cases: the first one is  
\(\sum_{s=1}^{a} i_s = \sum_{s=1}^{b} j_s = \sum_{s=1}^{c} k_s =1 \),  
the second one is  
 \(\{\sum_{s=1}^{a} i_s , \sum_{s=1}^{b} j_s , \sum_{s=1}^{c} k_s =1 \}
={0,1,2} \} \) as a set, and the last case is that one is distributed by 3 and
the others are distributed by 0.   
\subsection{3 in one of three series}
\subsubsection{3 in 0-series} About this part, we refer to the subsubsection
\ref{w=m:3}.

\subsubsection{3 in 1-series}
Suppose that 1-series has the number 3 or numbers 1,2. 
Then we have joined two tables. 

\kmcomment{
\begin{align*} 
  &  \begin{array}{c|c|c|c}
\text{\almostprimary} & \text{0-series} &\text{ 1-series} & \text{2-series } \\ \hline
\text{almost 2ndary} & \alpha & 3,0,\ldots,0 & \beta \\ \hline
\text{carnality} & \#\alpha & m - \#\alpha - \#\beta & \#\beta
\end{array} 
\\
  &   \begin{array}{c|c|c|c}
\text{\almostprimary} & \text{0-series} &\text{ 1-series} & \text{2-series } \\ \hline
\text{almost 2ndary} & 1,2, \alpha & 0,\ldots,0 & \beta \\ \hline
\text{carnality} & 2+ \#\alpha & m - 2- \#\alpha - \#\beta & \#\beta
\end{array} 
\end{align*} 
}

\begin{align*} 
    \begin{array}{c|c|c|c|*{3}{|c}}
        \hline
\text{\almostprimary} & \text{0-series} &\text{ 1-series} & \text{2-series }
                      & 
 \text{0-series} &\text{ 1-series} & \text{2-series } \\ \hline 
\text{almost 2ndary} &  \alpha & 3,0,\ldots,0 & \beta
                     & 
  \alpha & 1,2,0,\ldots,0 & \beta \\ \hline
\text{carnality} &  \#\alpha & m -  \#\alpha - \#\beta & \#\beta
 &  \#\alpha & m - \#\alpha - \#\beta & \#\beta \\\hline
\end{array} 
\end{align*}
We check the primary weight condition.  

Case 1. \( m-3 = w = m - \# \alpha + \# \beta \) and 
\(\# \alpha  = 3 + \# \beta)  \) and it is impossible. 

Case 2. \( m-3 = w = m - \# \alpha + \# \beta \) and 
\( \#\alpha  = 3 +  \#\beta\) and it is impossible. 

\subsubsection{3 in 2-series}
Suppose that 2-series has the number 3 or 1,2. Then we have joined two
tables.


\begin{align*} 
    \begin{array}{c|c|c|c|*{3}{|c}}
        \hline
\text{\almostprimary} & \text{0-series} &\text{ 1-series} & \text{2-series }
                      & 
 \text{0-series} &\text{ 1-series} & \text{2-series } \\ \hline 
\text{almost 2ndary} &  \alpha & 0,\ldots,0 & 3, \beta
                     & 
  \alpha & 0,\ldots,0 & 1,2, \beta \\ \hline
\text{carnality} &  \#\alpha & m - 1- \#\alpha - \#\beta & 1+ \#\beta
 &  \#\alpha & m - 2- \#\alpha - \#\beta & 2+ \#\beta \\\hline
\end{array} 
\end{align*}
We check the primary weight condition.  

Case 1. \( m-3 = w = m - 1 - \# \alpha + \# \beta\) and 
\( \#\alpha  = 2 +  \#\beta  \) and it is impossible. 

Case 2. \( m-3 = w = m +2 - \# \alpha + \# \beta\) and 
\(\# \alpha = 5 + \# \beta\) and it is impossible. 

\subsection{1 for the three series}
\begin{wrapfigure}[3]{r}[5mm]{0.55\textwidth}
\vspace{-1.0 \baselineskip}
\( \begin{array}{c|c|c|c}
\text{\almostprimary} & \text{0-series} &\text{ 1-series} & \text{2-series } \\ \hline
\text{almost 2ndary} & 1, \alpha & 1,\ldots,0 & 1, \beta \\ \hline
\text{carnality} & 1+ \#\alpha & m - 2- \#\alpha - \#\beta & \#\beta
\end{array} \)  \end{wrapfigure}
The primary weight becomes \(  
 m - 2- \#\alpha - \#\beta + 2+ 2 \#\beta
 = m + 2- \#\alpha+  \#\beta
 \) and we get 
 \( 0= 5- \#\alpha + \#\beta \) and no solution. 

\subsection{0,1,2 for all series} We let \(\{p,q\} =\{1,2\}\) as a set.  

\subsubsection{0 is in 0-series:}  
\begin{wrapfigure}[3]{r}[5mm]{0.55\textwidth}
\vspace{-1.5 \baselineskip}
\( \renewcommand{\arraystretch}{0.8}
\begin{array}{c|c|c|c}
\text{\almostprimary} & \text{0-series} &\text{ 1-series} & \text{2-series } \\ \hline
\text{almost 2ndary} & 0 & p,\ldots,0 & q, \beta \\ \hline
\text{carnality} & 1 & m - 2-  \#\beta & 1+  \#\beta
\end{array} \) 
\end{wrapfigure}
The primary weight becomes 
\(  m - 2- \#\alpha - \#\beta + 2+ 2 \#\beta
 = m  - \#\alpha + \#\beta
 \) and we get 
 \( 0= 2- \#\alpha + \#\beta \) and no solution. 

\subsubsection{0 is in 1-series:} 
\begin{wrapfigure}[3]{r}[5mm]{0.60\textwidth}
\vspace{-0.0 \baselineskip}
\( \renewcommand{\arraystretch}{0.8}
\begin{array}{c|c|c|c}
\text{\almostprimary} & \text{0-series} &\text{ 1-series} & \text{2-series } \\ \hline
\text{almost 2ndary} & p, \alpha & 0,\ldots,0 & q, \beta \\ \hline
\text{carnality} & 1+ \#\alpha & m - (1+ \#\alpha) - (1+ \#\beta) & 1+ \#\beta
\end{array} \) 
\end{wrapfigure}
The primary weight becomes \(  m - (1+ \#\alpha) + (1+ \#\beta) 
 = m  - \#\alpha + \#\beta
 \) and we get 
 \( 0= 3- \#\alpha + \#\beta \) and no solution. 

\subsubsection{ 0 is in 2-series:} 

\begin{wrapfigure}[3]{r}[5mm]{0.55\textwidth}
\vspace{-1.5 \baselineskip}
\( \begin{array}{c|c|c|c}
\text{\almostprimary} & \text{0-series} &\text{ 1-series} & \text{2-series } \\ \hline
\text{almost 2ndary} & p, \alpha & q,\ldots,0 & 0 \\ \hline
\text{carnality} & 1+ \#\alpha & m - (1+ \#\alpha) - 1  & 1
\end{array} \) 
\end{wrapfigure}
The primary weight becomes \(  
 m - (1+ \#\alpha) + 2 
 = m +1 - \#\alpha \) and we get 
 \( 0= 4- \#\alpha  \) and no solution. 

}



\begin{thebibliography}{1}

\bibitem{Mik:Miz:superForms}
Kentaro Mikami and Tadayoshi Mizutani.
\newblock Superalgebra structure on differential forms of manifold, 
\newblock{arXiv:2105.09738}, May 2021.

\bibitem{Mik:Miz:super3}
Kentaro Mikami and Tadayoshi Mizutani.
\newblock The second {B}etti number of doubly weighted homology groups of some
  pre {L}ie superalgebra.
\newblock {\em Tohoku Math. J.}, 74(2):301--311, 2022.

\end{thebibliography}

\newcommand{\noop}[1]{}\def\cprime{$'$} \def\cprime{$'$}

\end{document}